\documentclass[11pt,a4paper]{article}
\usepackage[T1]{fontenc}
\usepackage[utf8]{inputenc}
\usepackage{lmodern,microtype,amsmath,amssymb,amsthm,mathtools,booktabs,longtable,array,enumitem,xurl}
\usepackage[margin=2.55cm]{geometry}
\usepackage{xcolor}
\usepackage[colorlinks=true,linkcolor=blue!50!black,citecolor=blue!45!black,urlcolor=blue!55!black]{hyperref}
\setlist{leftmargin=*,itemsep=2pt,topsep=4pt}
\numberwithin{equation}{section}
\newtheorem{theorem}{Theorem}[section]
\newtheorem{lemma}[theorem]{Lemma}
\newtheorem{proposition}[theorem]{Proposition}
\newtheorem{corollary}[theorem]{Corollary}
\theoremstyle{remark}
\newtheorem{remark}[theorem]{Remark}
\newcommand{\F}{\mathbb F}
\newcommand{\Z}{\mathrm Z}
\hypersetup{pdftitle={The Exact Central-Index Profile of Groups with Bounded Noncommutativity},pdfauthor={Guillaume Lecomte}}
\title{\bfseries The Exact Central-Index Profile\\of Groups with Bounded Noncommutativity}
\author{Guillaume Lecomte\\[-1pt]
 \small Independent researcher\\[-1pt]
 \small \href{https://orcid.org/0009-0003-3532-3379}{ORCID: 0009-0003-3532-3379}\\[-1pt]
 \small \href{mailto:guillaume.lecomteexed@edu.executive.em-lyon.com}{guillaume.lecomteexed@edu.executive.em-lyon.com}}
\date{}
\begin{document}
\maketitle
\begin{abstract}
For a finite group $G$, let $\omega(G)$ be the largest size of a pairwise
noncommuting set and put $z(G)=[G:\Z(G)]$. We determine the maximum
possible central index under $\omega(G)\le k$ for every integer
$k\ge1$, giving an exact answer to Neumann's central-index question
and sharpening Pyber's exponential upper bound.

The maximum is $1$ for $k=1,2$, $(k-1)^2$ for $3\le k\le6$,
$81$ for $k=7,8$, and $4^{\lfloor(k-1)/2\rfloor}$ for $k\ge9$.
In the last range, with $r=\lfloor(k-1)/2\rfloor$, equality holds
exactly when $G'\cong C_2$ and $G/\Z(G)\cong\F_2^{2r}$.
We also obtain the sharp universal bound
$z(G)\le(81/64)2^{\omega(G)-1}$.

All maximizers are classified up to isoclinism. The results extend
to arbitrary groups with finite clique number through finite
isoclinic representatives.
\end{abstract}

\section{Introduction}\label{sec:intro}
The size of a pairwise noncommuting set measures how many elements
of a group can be mutually incompatible with commutativity. For a
finite group $G$, let $\omega(G)$ be the largest size of such a set.
We ask how strongly a bound on $\omega(G)$ constrains the index
of the centre. The order of $G$ itself cannot be bounded in this
way: taking a direct product with an arbitrarily large abelian
group changes neither $\omega(G)$ nor the central index. The
natural extremal quantities are therefore
\[
 z(G)=[G:\Z(G)],\qquad
 \zeta(k)=\sup\{z(G):\omega(G)\le k\}.
\]
Here the supremum is over finite groups. In graph-theoretic
language, join two elements of $G$ when they do not commute;
then $\omega(G)$ is the clique number of this noncommuting graph.
Thus the problem relates a combinatorial restriction to a
structural invariant of the group.

The origin is a question posed by Erd\H{o}s in 1975: if a group
has no infinite pairwise noncommuting subset, must the sizes of
its finite pairwise noncommuting subsets be bounded?
Neumann answered this affirmatively in 1976 by proving the
stronger statement that such a group has finite central quotient
\cite[Theorem 6]{Neumann1976}. We first solve the finite extremal
problem and then extend the sharp bounds to arbitrary groups by
isoclinism, as explained in Corollary~\ref{cor:arbitrary} below.
Neumann also explicitly asked how the central index can be bounded
in terms of the largest noncommuting set. His proof gave
$\log z(G)=O(\omega(G)^2)$, and he suggested a possible quadratic
bound for $z(G)$ itself \cite[Section 4]{Neumann1976}.

Pyber resolved the scale of this quantitative question in 1987:
there is an absolute constant $c$ such that
$z(G)\le c^{\omega(G)}$ for every finite group
\cite{Pyber1987}. In particular, $\zeta(k)$ is finite.
The extraspecial $2$-groups show that exponential growth is
necessary; their precise parameters are recalled below.
The remaining extremal questions are to determine the optimal
exponential base, the exact maximum at each finite cap, and the
groups attaining it. These are distinct from the related
questions about large abelian subgroups and the number of
commuting tuples studied by Erd\H{o}s and Straus
\cite{ErdosStraus1976}.

Several parts of the literature supply information about the
models and local estimates relevant to these questions.
Isaacs's formula, reported by Bertram, gives the clique numbers of
extraspecial $2$-groups \cite[p.~40]{Bertram1983}; Chin studies noncommuting sets in
extraspecial $p$-groups \cite{Chin2005}. Further results treat
groups with derived subgroup of prime order and metacyclic
$p$-groups \cite{LiuWang2015,FouladiOrfi2013}.
For scalar commutator maps, noncommuting sets correspond to sets
of pairwise nonorthogonal points in symplectic geometry. Results
on partial ovoids therefore provide sharp information in some
small odd-characteristic cases
\cite{BrouwerVanMaldeghem2022,Ceria2023}.
Covering theorems provide a complementary tool: the centralizers
of the elements of a maximum clique cover the group, so bounds
for subgroup and coset covers can control the central index
\cite{Neumann1954Cosets,Tomkinson1987,BryceFedriSerena1997}.

The present article determines $\zeta(k)$ and classifies all its
maximizers up to isoclinism for every $k\ge1$. The eventual
extremizers are precisely the extraspecial binary isoclinism
classes, but small parameters also require quinary and ternary
models. The proof combines the geometry of the nilpotent cases
with a centralizer descent that excludes larger values in
arbitrary finite groups. We first state the exact profile, then
describe its models, its extension to arbitrary groups, and the
structure of the proof.

\begin{theorem}[The central-index profile]\label{thm:global}
For every integer $k\ge1$, the supremum defining $\zeta(k)$ is
attained and equals
\begin{equation}\label{eq:profile}
 M(k)=
 \begin{cases}
 1,&k=1,2,\\
 (k-1)^2,&3\le k\le6,\\
 81,&k=7,8,\\
 4^{\lfloor(k-1)/2\rfloor},&k\ge9.
 \end{cases}
\end{equation}
\end{theorem}

The three nontrivial ranges have different extremal models. In the
Heisenberg group $\mathrm{UT}_3(\F_q)$, the commutator on the
central quotient $\F_q^2$ is the determinant. Pairwise
noncommuting elements project to distinct lines, and one may choose
one representative of every line. Consequently
\[
 z\bigl(\mathrm{UT}_3(\F_q)\bigr)=q^2,
 \qquad \omega\bigl(\mathrm{UT}_3(\F_q)\bigr)=q+1.
\]
Taking $q=2,3,4,5$ explains the four initial squares. An extraspecial
3-group of order $3^5$ instead has $z=81$ and $\omega=7$.
Its cliques are partial ovoids of $W(3,3)$, that is, sets of
pairwise nonorthogonal projective points: Tallini's bound
$q^2-q+1$ for odd $q$ gives the upper bound seven, and two
symplectic planes give the matching construction
\cite[Section 2 and Table 1]{Ceria2023}; see also
\cite{BrouwerVanMaldeghem2022,CimrakovaFack2005}.
Finally, for an extraspecial $2$-group $E_r$ of order $2^{1+2r}$,
Isaacs's clique formula, reported by Bertram
\cite[p.~40]{Bertram1983}, gives
\[
 z(E_r)=4^r,\qquad \omega(E_r)=2r+1.
\]
Thus the largest allowed $r$ gives
$L_k=4^{\lfloor(k-1)/2\rfloor}$ and explains the paired plateaux.
These examples prove all lower bounds in Theorem~\ref{thm:global}.

\begin{theorem}[All extremizers]
\label{thm:equality}
Let $G$ be finite with $\omega(G)\le k$.
\begin{enumerate}[label=\textup{(\roman*)}]
\item If $k=1$ or $2$, then $G$ is abelian and $z(G)=1$.
\item If $k\in\{3,4,6\}$ and $p=k-1$, then $z(G)=p^2$ if and only if
$G'\cong C_p$ and $G/\Z(G)\cong\F_p^2$; such groups have
$\omega(G)=k$ and are isoclinic to an extraspecial group of order $p^3$.
\item If $k=5$, then $z(G)=16$ if and only if $G$ is isoclinic
to either an extraspecial group of order $32$ or
$\mathrm{UT}_3(\F_4)$. These are two distinct isoclinism classes,
and both have $\omega(G)=5$.
\item If $k=7$ or $8$, then $z(G)=81$ if and only if
$G'\cong C_3$ and $G/\Z(G)\cong\F_3^4$; such groups have
$\omega(G)=7$.
\item If $k\ge9$ and $r=\lfloor(k-1)/2\rfloor$, then
$z(G)=L_k$ if and only if
$G'\cong C_2$ and $G/\Z(G)\cong\F_2^{2r}$; such groups have
$\omega(G)=2r+1$.
\end{enumerate}
In \textup{(iv)} and \textup{(v)}, the groups form the corresponding
extraspecial isoclinism class. Thus all maximizers are identified
up to isoclinism for every cap $k$.
\end{theorem}

\begin{corollary}[Sharp exponential bound]\label{cor:universal}
Every finite group satisfies
\[
 z(G)\le\frac{81}{64}\,2^{\omega(G)-1}.
\]
Equality holds precisely when $G'\cong C_3$ and
$G/\Z(G)\cong\F_3^4$. Moreover,
\begin{equation}\label{eq:exponential-consequences}
 z(G)\le2^{\omega(G)-1}\quad\text{if }\omega(G)\ge9,
 \qquad \lim_{k\to\infty}\zeta(k)^{1/k}=2.
\end{equation}
\end{corollary}

\begin{corollary}[Arbitrary groups]\label{cor:arbitrary}
Let $G$ be a group, not necessarily finite, with $\omega(G)\le k<\infty$.
Then $[G:\Z(G)]\le M(k)$, and all the equality characterizations
in Theorem~\ref{thm:equality} remain valid.
The sharp bound and its equality case in Corollary~\ref{cor:universal}
also remain valid when $\omega(G)<\infty$.
For every $k$, the bound $M(k)$ is attained by both finite and
infinite groups.
\end{corollary}
\begin{proof}
The abelian case is immediate. Otherwise, Zarrin's finite-representative
theorem \cite[Theorem 1.1\textup{(2)}]{Zarrin2016} gives a finite
group $H$ isoclinic to $G$, with $\omega(H)=\omega(G)$.
Isoclinism also gives $H/\Z(H)\cong G/\Z(G)$ and $H'\cong G'$.
Thus $z(H)=z(G)$, and the finite results give every asserted bound
and equality characterization. Taking the direct product of a finite
maximizer with an infinite abelian group preserves both invariants
and proves the last assertion.
\end{proof}

Underlying this established reduction are the finiteness of the
central quotient, Schur's theorem, and Hall's existence of a stem
representative $H$, meaning that $\Z(H)\le H'$ \cite{Hall1940}.
Isoclinism preserves the commutator map and hence the clique number
\cite[Lemma 2.1 and proof of Theorem 1.1]{Zarrin2016}.
In particular, the equality statements specify an isoclinism class;
they do not assert a direct-product decomposition of an infinite
extremizer.

\paragraph{Proof strategy.}
The proof separates the nilpotent structure from the passage to a
general group. Section~\ref{sec:clique-calculus} develops one
calculus for both: cliques in cosets and commuting products.
Section~\ref{sec:scalar-budgets} proves the Sylow budgets, the
complete nilpotent profile and the binary rigidity theorem.

For a centralizer $K=C_G(L)$, the accounting identity is
\begin{equation}\label{eq:intro-accounting}
 z(G)=[G:K]\,z(K)\,[\Z(K):\Z(G)].
\end{equation}
Section~\ref{sec:descent} introduces the prime cost of losing
centrality and then uses Baer--Suzuki to choose a
nonnilpotent $L$, whose central quotient amplifies cliques of $K$
four times. Two saturation principles close the finite exceptions:
near-minimal ternary losses force divisibility, and affine
conjugacy classes force a faithful action on the class.

Section~\ref{sec:completion} first settles the finite caps and
then begins one induction at eighteen. A numerical prime-cost
bound pays for the internal centralizer index with the same
budget as the central loss. This leaves one centralizer estimate,
a finite list of saturation obstructions and one envelope for all
larger caps.
The extremizer statements then follow from the nilpotent rigidity
already proved, and the sharp exponential estimate follows from
the exact profile.

\paragraph{Notation.}
Groups in the proof of the finite results are finite;
Corollary~\ref{cor:arbitrary} explicitly removes that hypothesis.
For an arbitrary group, $\omega(G)$ denotes the supremum of the
cardinalities of pairwise noncommuting subsets; whenever it is
finite, it is a maximum. We write $G'$ for the derived subgroup
and $C_G(X)$ for the pointwise centralizer of a subset $X$, with
$C_G(u,v)=C_G(u)\cap C_G(v)$. We use $[x,y]=x^{-1}y^{-1}xy$ and
$x^g=g^{-1}xg$. A \emph{clique} is a pairwise noncommuting set;
\emph{maximum} means of largest cardinality. We also use $\omega(X)$
for a subset or a coset. Two sets are \emph{completely adjacent}
if every element of one fails to commute with every element of
the other. An abelian group has clique number one,
and no group has clique number two: a noncommuting pair $x,y$
extends to the clique $x,y,xy$. The integer $k\ge1$ is an upper
cap; for a fixed ambient group $G$, write $N=\omega(G)$ and
$\Z=\Z(G)$. Logarithms are binary unless specified.
$C_m$ denotes the cyclic group of order $m$ and $\F_q$ the field
with $q$ elements. We write $W(2r-1,q)$ for the symplectic polar
space associated with $\F_q^{2r}$. For a finite abelian or nilpotent
group $H$, the symbol $H_{(p)}$ denotes its Sylow $p$-subgroup.
We use $\gamma_i(P)$ for the lower central series and
$\kappa=81/64$ for the sharp constant above. Isoclinism means an
identification of the central quotients and derived groups
compatible with their commutator maps; it does not require a
direct-product decomposition with an extraspecial factor.

\section{Clique constructions}\label{sec:clique-calculus}
We use three ways to construct cliques: inside a coset, by adjoining
an element that moves a central subgroup, and across commuting
factors. Only the coset estimate requires a covering argument.

\subsection{Cliques in a coset}
A coset cover is \emph{irredundant} if none of its members can be
deleted. A maximum clique in a coset will produce such a cover by
centralizers; the next estimate turns that cover into a lower
bound for the clique.
For its context see \cite{Neumann1954Cosets,Tomkinson1987,
BryceFedriSerena1997,Berkovich2010}; the binary vector-space case
is also covered by \cite{Brouwer1986}.

\begin{lemma}[Coset-covering bound for a $p$-group]
\label{lem:coset-cover}
Let $U$ be a finite $p$-group and
$U=\bigcup_{i=1}^t H_i x_i$ an irredundant covering by right
cosets of subgroups. Then
\[
 [U:\bigcap_{i=1}^t H_i]\le p^{t-1}.
\]
\end{lemma}
\begin{proof}
We induct on $|U|$, the case $t=1$ being immediate. For $t>1$
all $H_i$ are proper. Choose a maximal subgroup $M\ge H_1$;
it is normal of index $p$. Translate the cover so that
$H_1x_1\cap M=\varnothing$, and choose an irredundant subcover
of $M$ with indices $I_1\subseteq\{2,\ldots,t\}$. Put
$B_1=\bigcap_{i\in I_1}H_i$ and $D=\bigcap_iH_i$.
Induction gives $[M:M\cap B_1]\le p^{|I_1|-1}$.

Suppose the subgroup intersection for the selected indices is $B_j\ne D$.
Their cosets form a proper subfamily of the original irredundant
cover, so choose $y$ outside their union. The coset $B_jy$
avoids them all. Intersect the remaining cosets with $B_jy$
and translate by $y^{-1}$; their nonempty intersections are
right cosets in $B_j$, with directions $B_j\cap H_i$.
An irredundant subcover selects new indices $I_{j+1}$ and gives
\[
 B_{j+1}=B_j\cap\bigcap_{i\in I_{j+1}}H_i,
 \qquad [B_j:B_{j+1}]\le p^{|I_{j+1}|-1}.
\]
Induction applies since $B_j\le B_1<U$. The nonempty index sets
are disjoint, so this process reaches $B_s=D$ after finitely
many steps, even if some successive subgroups coincide.
Since $D\le H_1\le M$,
\[
 [U:D]\le p[M:M\cap B_1][B_1:D]
 \le p^{1+\sum_{j=1}^s(|I_j|-1)}\le p^{t-1}.
\]
For the last inequality, $s=1$ gives $|I_1|\le t-1$;
if $s\ge2$, use $\sum_j|I_j|\le t$.
No subgroup direction is assumed normal.
\end{proof}

\begin{lemma}[Coset calculus]\label{lem:coset-calculus}
Let $U\le G$, $g\in G$, and let $X$ be a maximum clique in $gU$.
Then
\[
 U\cap C_G(X)=\Z(U)\cap C_U(g).
\]
The sets $C_G(x)\cap gU$, $x\in X$, form an irredundant coset
cover, and some $x\in X$ satisfies
$[U:C_U(x)]\le|X|$. If $U$ is a $p$-group, then
\begin{equation}\label{eq:whole-coset}
 \omega(gU)\ge1+\log_p[U:\Z(U)\cap C_U(g)].
\end{equation}
If $U$ is abelian, the exact formula is
\begin{equation}\label{eq:abelian-coset}
 \omega(gU)=[U:C_U(g)].
\end{equation}
If $U$ is a nonabelian $p$-group, then $\omega(gU)\ge p+1$.
\end{lemma}
\begin{proof}
The right-hand side of the first identity centralizes all of $gU$.
Conversely, suppose $a\in U$ centralizes $X$ but not some $y\in gU$.
Multiply by $a$ precisely the vertices of $X$ commuting with $y$.
Mutual noncommutation is preserved, and the new clique extends by
$y$, a contradiction. This proves the identity.

Maximality makes the displayed sets cover $gU$, and each $x\in X$
belongs only to its own set, so the cover is irredundant.
The part at $x$ is $xC_U(x)$ and has size $|C_U(x)|$;
counting gives the averaging assertion. Translation by $g^{-1}$
identifies the parts with left cosets in $U$. Inversion turns them
into right cosets with the same subgroup directions. The intersection
of those directions is $\Z(U)\cap C_U(g)$, so
Lemma~\ref{lem:coset-cover} proves \eqref{eq:whole-coset}.
Finally, when $U$ is abelian, $ga$ and $gb$ commute exactly when
$ab^{-1}\in C_U(g)$, proving \eqref{eq:abelian-coset}.

For a nonabelian $p$-group, a one-part cover would force
$U=\Z(U)\cap C_U(g)$. Otherwise every direction is proper and
each part has size at most $|U|/p$.
At most $p$ parts could cover only by partitioning into $p$
index-$p$ cosets. Their directions are normal, and cosets of
distinct index-$p$ subgroups intersect, so all directions coincide.
Their intersection lies in $\Z(U)$, making $U/\Z(U)$ cyclic,
a contradiction.
\end{proof}

\begin{lemma}[Clique extensions]\label{lem:twist}
If $A\le\Z(H)$ and $[A,g]\ne1$, then
\begin{equation}\label{eq:twisting}
 \omega(\langle H,g\rangle)\ge\omega(H)+[A:C_A(g)].
\end{equation}
If $[K,L]=1$, then
\begin{equation}\label{eq:product-amplification}
 \omega(KL)\ge\omega(K)\omega(L/\Z(L)).
\end{equation}
\end{lemma}
\begin{proof}
The abelian-coset formula gives $[A:C_A(g)]$ vertices in $gA$.
Choose $a_0\in A\setminus C_A(g)$. In a maximum clique of $H$,
multiply by $a_0$ precisely the vertices commuting with $g$.
Their mutual commutators are preserved, and every resulting vertex
is noncommuting with every element of $gA$. This proves
\eqref{eq:twisting}.

For \eqref{eq:product-amplification}, multiply a maximum clique
of $K$ by lifts of a clique in $L/\Z(L)$. Cross commutators from
$L$ lie outside $\Z(L)$, whereas $K\cap L\le\Z(L)$, so they
cannot cancel commutators from $K$. Products using the same lift
remain noncommuting inside $K$.
\end{proof}

\begin{lemma}[Small central quotients]\label{lem:small-quotients}
If $Q$ is nonabelian and $\omega(Q)\le4$, then
$Q/\Z(Q)$ is $C_2^2$, $C_3^2$, or $S_3$. Thus $z(Q)=4$ for
$\omega(Q)=3$, and $z(Q)\le9$ for $\omega(Q)=4$.
\end{lemma}
\begin{proof}
For $x\notin\Z(Q)$, twisting on $\langle x\rangle\le\Z(C_Q(x))$
gives $\omega(C_Q(x))\le\omega(Q)-2\le2$. Every proper centralizer is
therefore abelian. Their noncentral parts partition $Q\setminus\Z(Q)$,
and their number is $s=\omega(Q)\in\{3,4\}$. Write their orders modulo
$\Z(Q)$ as $m_1\ge\cdots\ge m_s\ge2$ and put $q=z(Q)$. Then
\[
 q=1+\sum_i(m_i-1),\qquad
 m_1m_2\le q\le m_1+(s-1)(m_2-1),
\]
so $m_1\le s-1$. For $s=3$ all orders are two, giving $q=4$.
For $s=4$, the conditions $m_i\mid q$ leave only $(3,2,2,2)$ or
$(3,3,3,3)$, giving $q=6$ or nine. A cyclic central quotient is
impossible; the listed isomorphism types follow.
\end{proof}

\section{The nilpotent profile and rigidity}\label{sec:scalar-budgets}
The common coset calculus controls every Sylow subgroup through a
scalar commutator form and the cliques in its fibres. We prove the
prime budgets first, combine them into the nilpotent profile, and
then establish the binary gap that will classify the maximizers.
We use the elementary fact that if $G/\Z(G)$ is a $p$-group,
lifting its central series makes $G$ nilpotent and all Sylow
subgroups at other primes central.

\subsection{The symplectic model}
For a symplectic form $b$, a clique means a set whose distinct vectors
have nonzero $b$-pairing. Such sets model noncommuting elements
\cite{Bertram1983,Chin2005,LiuWang2015}. The hyperplane variant below
will also be used in the equality proof.

\begin{lemma}[Symplectic cliques]\label{lem:symplectic-cliques}
Let $(V,b)$ be symplectic of dimension $2r>0$ over $\F_p$.
There is a clique of size $pr+1$, and one can choose such a clique
to meet any prescribed hyperplane in exactly one point.
Every maximum clique spans $V$.
For $p=2$, the maximum is exactly $2r+1$, and every maximum
clique has sum zero.
For $p=3$, a clique of size at least $4r-1$ exists.
\end{lemma}
\begin{proof}
In a symplectic plane, choose one representative of each of its
$p+1$ lines. If $X,Y$ are cliques in orthogonal summands and
$y_0\in Y$, then $(X+y_0)\cup(Y\setminus\{y_0\})$ is a
clique of size $|X|+|Y|-1$. This proves the lower bound $pr+1$.
For the hyperplane assertion write $J=e^\perp$, choose
$b(e,f)=1$, and split $V=\langle e,f\rangle\perp V_0$.
For a clique $D_0\subset V_0$ of size $p(r-1)+1$, using
$D_0=\{0\}$ when $r=1$, the set
\[
 (D_0+f)\cup\{e+tf:t\in\F_p\}
\]
is a clique of size $pr+1$, and only $e$ belongs to $J$.

If a maximum clique $D$ does not span, choose
$0\ne v\in\langle D\rangle^\perp$ and $y$ with $b(v,y)\ne0$.
Translate by $v$ precisely the vertices of $D$ orthogonal to $y$.
All mutual pairings are unchanged, and the resulting clique extends
by $y$, a contradiction.

Over $\F_2$, a clique of size $t$ has Gram matrix $J_t-I_t$,
where $J_t$ has all entries one and $I_t$ is the identity matrix.
The Gram matrix has rank $t$ for even $t$ and $t-1$ for odd $t$.
Hence $t\le2r+1$, with equality by the construction above.
The sum of a maximum clique is orthogonal to each of its vertices;
as these span, the sum is zero.
For $p=3$, \cite[Theorem~3.7]{Ceria2023} supplies a partial ovoid
of $W(5,3)$ of size $3^2+3+1=13$. Choosing one nonzero vector
on each projective point gives thirteen vertices in dimension six.
For $r=3u+v$, $0\le v\le2$,
orthogonal gluing with $v$ planes gives
$12u+3v+1\ge4r-1$ vertices.
\end{proof}

The next lemma enlarges cliques when a scalar form hides a second
commutator coordinate. It will be used in the small ternary case
and in the binary rigidity proof.

\begin{lemma}[Affine inflation of a nonscalar form]
\label{lem:affine-inflation}
Let $\beta:V\times V\to W$ be alternating bilinear over $\F_p$,
with values spanning $W$, $\dim V=2r$ and $\dim W\ge2$.
If some scalarization $b=\lambda\circ\beta$ is nondegenerate, then
$r\ge2$ and
\[
 \omega(\beta)\ge p^2(r-1)+1,
\]
where $\omega(\beta)$ is the largest size of a set with pairwise
nonzero $\beta$-values.
\end{lemma}
\begin{proof}
Choose $\mu\in W^*$ independent of $\lambda$ and put
$c=\mu\circ\beta$. Define $T$ by $c(v,x)=b(Tv,x)$.
Since the values of $\beta$ span $W$, $T$ is not scalar.
Hence some $v,Tv$ are linearly independent: an operator preserving
every line is scalar, as a basis and its pairwise sums show.
On $U=v^{\perp_b}$ the functional $x\mapsto c(v,x)$ is therefore
nonzero and descends to $S=U/\langle v\rangle$.
This space is symplectic of dimension $2r-2$, so $r\ge2$.

By Lemma~\ref{lem:symplectic-cliques}, $S$ has a scalar clique
of size $p(r-1)+1$ meeting the kernel hyperplane once.
Lift its vertices to $x_0,x_1,\ldots,x_{p(r-1)}\in U$,
with $x_0$ exceptional. Replace each $x_i$, $i>0$, by the
affine line $x_i+\F_pv$. Within a line, distinct vectors have
nonzero $\beta$-pairing because $c(v,x_i)\ne0$.
Between different lines, and with $x_0$, the $b$-pairings
remain the nonzero quotient pairings. The resulting clique
has $1+p\cdot p(r-1)$ vertices.
\end{proof}

\subsection{Scalar fibres and Sylow budgets}
Scalarization applies a homomorphism with values in $\F_p$ to the
commutator. A nonzero scalar value certifies noncommutation, but a
zero value need not certify commutation. The radical consists of
those elements invisible to the chosen scalar form; we must therefore
count cliques inside its cosets as well as between them.

Let $P$ be a nonabelian $p$-group. A nonzero homomorphism
$\chi:P'/\gamma_3(P)\to\F_p$ exists by nilpotence. We also write
$\chi$ for its composite with $P'\to P'/\gamma_3(P)$. The pairing
\[
 b_\chi(x,y)=\chi([x,y])
\]
is biadditive and conjugation-invariant. Its radical is
\[
 R=R_\chi=\{x\in P:b_\chi(x,y)=0\ \text{for all }y\in P\}.
\]
It is normal and contains $P'P^p\Z(P)$. The induced form on
$V=P/R$ is nondegenerate, making $V$ an elementary abelian
symplectic space. All scalarizations below use a nonzero $\chi$.
We use the parameters
\begin{equation}\label{eq:scalar-parameters}
 \begin{gathered}
 d=\dim_{\F_p}V=2r,\qquad A=\Z(R),\qquad
 m=\log_p[R:A],\qquad s=\log_p[A:\Z(P)],\\
 n=\log_pz(P)=d+m+s,\qquad h=\log_p|[A,P]|.
 \end{gathered}
\end{equation}
Since $R$ is normal in $P$, its centre $A=\Z(R)$ is normal and abelian.
The identity $n=d+m+s$ separates the scalar quotient from two
parts hidden in its radical: noncommutativity inside $R$, measured
by $m$, and the excess of $\Z(R)$ over $\Z(P)$, measured by $s$.
The term $h$ records the action of $P$ on $\Z(R)$.
Here $m=0$ or $m\ge2$, since a nontrivial central quotient is
not cyclic, and $[A,P]$ is the subgroup generated by the commutators.
We call a coset $gR$ a \emph{scalar fibre}. Fibres above distinct
vertices of a scalar clique are completely noncommuting with one
another, so their internal clique sizes can be added.

The odd-prime estimates apply in every nilpotency class; compare
the restricted families in
\cite{AzadFouladiOrfi2013,FouladiOrfi2013,LiuWang2015}.

\begin{theorem}[Sylow budgets]\label{thm:sylow-budgets}
Let $P$ be a finite $p$-group, and put $n=\log_pz(P)$.
For nonabelian $P$, use the scalarization parameters in
\eqref{eq:scalar-parameters}.
\begin{enumerate}[label=\textup{(\roman*)}]
\item If $p=2$, every scalarization of nonabelian $P$ satisfies
\begin{equation}\label{eq:binary-excess}
 \omega(P)\ge n+1+dm+h.
\end{equation}
For nonabelian $P$, equality $\omega(P)=n+1$ forces $P$ to have
class two and $P/\Z(P)$ to be elementary abelian of even dimension.
Every $2$-group consequently satisfies
\begin{equation}\label{eq:binary-exact}
 z(P)\le L_{\omega(P)}.
\end{equation}
\item For odd $p$,
\begin{equation}\label{eq:odd-budgets}
 \omega(P)-1\ge\begin{cases}pn/2,&p\ge5,\\2n-2,&p=3.\end{cases}
\end{equation}
In particular, $z(P)\le2^{\omega(P)-1}$ for $p\ge5$, and
$z(P)\le\kappa2^{\omega(P)-1}$ for $p=3$.
\item If $p=3$ and $z(P)=27$, then $\omega(P)\ge10$.
If $p=3$, $z(P)=81$ and $\omega(P)\le8$, then
$P'\cong C_3$, $P/\Z(P)\cong\F_3^4$ and $\omega(P)=7$.
\end{enumerate}
\end{theorem}
\begin{proof}
The abelian case is immediate, so assume $P$ nonabelian.
The common fibre estimate below is the input for all three parts.
The binary case uses a relation between the conjugation actions;
the odd cases use lower bounds for the sizes of the fibres.

\emph{Common fibre estimate.}
Choose lifts $g_v$ of a spanning scalar clique $D\subset P/R$
and put $\delta_v=\log_p[A:C_A(g_v)]$.
The coset bounds give
$\omega(g_vR)\ge m+\delta_v+1$ and
$\omega(g_vR)\ge\omega(g_vA)=p^{\delta_v}$.
Distinct fibres above $D$ are completely adjacent.
Since $A$ is normal and abelian, $[a,g_v]\in A$ and
$[aa',g_v]=[a,g_v][a',g_v]$ for $a,a'\in A$.
Thus the maps $a\mapsto[a,g_v]$ are homomorphisms, with common kernel
$\Z(P)$, because the lifts and $R$ generate $P$ and $R$ centralizes $A$.
Their joint image has order $p^s$ and embeds in the product of their
images, of orders $p^{\delta_v}$. Hence
\begin{equation}\label{eq:p-fibre-budget}
 \sum_{v\in D}\delta_v\ge s,\qquad
 \omega(P)\ge\sum_{v\in D}
 \max\{m+\delta_v+1,p^{\delta_v}\}.
\end{equation}

\emph{Part (i): binary interaction.}
For $p=2$, choose a zero-sum scalar clique of size $d+1$ and lift it to
$g_1,\ldots,g_{d+1}$; put $\delta_i=\log_2[A:C_A(g_i)]$.
Their actions $\alpha_i$ on $A$ commute,
since $R$ centralizes $A$, and $\prod_i\alpha_i=1$ because
$\prod_i g_i\in R$. Write $A$ additively and put
$A_i=(\alpha_i-1)A$ and $\beta_i=\alpha_1\cdots\alpha_i$, with
$\beta_0=1$.
The maps
\[
 \Phi:A\longrightarrow\bigoplus_i A_i,
 \quad a\longmapsto((\alpha_i-1)a)_i,
 \qquad
 T:\bigoplus_i A_i\longrightarrow A,
 \quad (v_i)\longmapsto\sum_i\beta_{i-1}v_i
\]
satisfy $\ker\Phi=\Z(P)$ and $T\Phi=0$ by telescoping.
Each $A_i$ is invariant under all $\beta_j$, so
$\operatorname{im}T=\sum_iA_i=[A,P]$.
Counting $\operatorname{im}\Phi\le\ker T$ gives
\begin{equation}\label{eq:telescope}
 \sum_i\delta_i\ge s+h.
\end{equation}
Substitution in \eqref{eq:p-fibre-budget} yields
$\omega(P)\ge(d+1)(m+1)+s+h=n+1+dm+h$.
If $\omega(P)=n+1$, then $m=h=0$, so $R=A=\Z(P)$.
Thus $P/\Z(P)$ is elementary abelian of even dimension and
$P'\le\Z(P)$. In general $n\le\omega(P)-1$; when $n$ is odd,
equality is impossible, so $n\le\omega(P)-2$.
This proves \eqref{eq:binary-exact}.

\emph{Part (ii): odd-prime budgets.}
Let $p$ be odd and set $\mu=m+1$, so $\mu=1$ or $\mu\ge3$.
Choose $D$ maximum, hence spanning. Its size is at least $pr+1$;
for $p=3$ it is at least $4r-1$, and is four when $r=1$.
For integers $\delta\ge0$, put
$f_p(\mu,\delta)=\max\{\mu+\delta,p^\delta\}$. Then
\[
 f_p(\mu,\delta)\ge\frac{\mu+1}{2}+\frac p2\delta,
 \qquad
 f_3(\mu,\delta)\ge\frac{\mu+1}{2}+2\delta
 \quad(\mu=1\text{ or }\mu\ge3).
\]
The first inequality averages the two arguments of the maximum,
using $p^\delta\ge1+(p-1)\delta$.
For the second, check $\delta=0,1$ directly; for $\delta\ge2$,
use the same averaging with $3^\delta\ge1+3\delta$.
Together with \eqref{eq:p-fibre-budget}, these give, for $p\ge5$,
\[
 \omega(P)\ge1+\frac{pn}{2}
 +\frac{(\mu-1)(p(r-1)+1)}2.
\]
For $p=3$, they give
\[
 \omega(P)\ge
 \begin{cases}
 2n,&r=1,\\[2pt]
 2n-1+\dfrac{(\mu-1)(4r-5)}2,&r\ge2.
 \end{cases}
\]
These prove \eqref{eq:odd-budgets}.
For $p\ge5$, $p/2>\log_2p$ gives $z(P)\le2^{\omega(P)-1}$.
For $p=3$ and $n\ge4$,
\[
 \frac{z(P)}{2^{\omega(P)-1}}\le4\left(\frac34\right)^n\le\frac{81}{64}.
\]
When $n=2$, the scalar-plane clique gives $\omega(P)\ge4$, hence ratio
at most $9/8$. When $n=3$, necessarily $r=1$, so
$\omega(P)\ge6$ and the ratio is at most $27/32$.
The case $n=1$ is impossible because a nontrivial central quotient
cannot be cyclic. This completes~\textup{(ii)}.

\emph{Part (iii): the small ternary gap.}
Suppose $p=3$, $[P:R]=9$ and $R\ne\Z(P)$.
Choose one nonzero representative of each of the four lines of $P/R$.
If $R$ is nonabelian, the nonabelian-coset bound gives four vertices
in each of the four completely adjacent fibres, hence sixteen.
If $R$ is abelian, the action of $P/R$ on $R$ is nontrivial.
Its kernel is a proper subspace containing at most one of our four
representatives. Each of the other three cosets has a clique of
size at least three; the fourth has at least one. Thus
$\omega(P)\ge10$ in either case.

At central index 27 every nonzero scalarization has rank two and
a noncentral radical, proving the first assertion of~\textup{(iii)}.
At index 81 with $\omega(P)\le8$, the same argument excludes
scalar radicals of index nine. Hence every scalar radical is the
centre, and $P$ has class two with $P/\Z(P)\cong\F_3^4$.
Also $P'$ is elementary abelian, since $[x,y]^3=[x^3,y]=1$.
If its dimension were at least two,
Lemma~\ref{lem:affine-inflation}, with $p=3$ and $r=2$, would give
ten vertices. Thus the commutator is a nondegenerate scalar form;
its clique number is seven by the $W(3,3)$ calculation recalled
in Section~\ref{sec:intro}. This proves~\textup{(iii)}.
\end{proof}

\subsection{Combining the Sylow budgets}
\begin{proposition}[The full nilpotent profile]\label{prop:nilpotent-profile}
If $F$ is nilpotent and $\omega(F)\le k$, then $z(F)\le M(k)$.
Every value $M(k)$ is attained by a nilpotent group.
For $k\ge9$, the comparison with $L_k$ is strict if an odd-order
Sylow subgroup of $F$ is nonabelian.
\end{proposition}
\begin{proof}
The abelian case is immediate. Write $F=\prod_pF_{(p)}$ and put
$t_p=\omega(F_{(p)})$. Cartesian products of cliques give
$\omega(F)\ge\prod_p t_p$. The Sylow budgets give
\[
 z(F)\le\kappa\,2^{\sum_p(t_p-1)},
\]
where the factor $\kappa$ is needed only in the ternary part.

First let $k=2a+1$ or $2a+2$, with $a\ge4$.
If exactly one Sylow subgroup is nonabelian, put
$n=\log_pz(F)$ for its prime $p$.
The binary bound gives $z(F)\le L_k$ when $p=2$.
For $p=3$, the odd budget gives $n\le a+1$, and
$3^{a+1}<4^a$. For $p\ge7$,
\[
 \log_2z(F)\le\frac{(4a+2)\log_2p}{p}
 <\frac{12a+6}{7}\le2a.
\]
For $p=5$, $5^3<2^7$ gives
$\log_2z(F)<(28a+14)/15\le2a$ for $a\ge7$;
for $a=4,5,6$, the respective bounds
$z(F)\le125,625,3125<4^a$ finish the comparison.
If at least two Sylow subgroups are nonabelian, their clique
numbers are at least three and one is at least four. Hence
\[
 \prod_p t_p-\sum_p(t_p-1)\ge7,
 \qquad z(F)\le\kappa2^{\omega(F)-7}<2^{k-2}\le L_k.
\]

The small quotients in Lemma~\ref{lem:small-quotients} handle
$k\le4$. For $k=5,6,8$, at most one Sylow subgroup is nonabelian,
since two would already supply at least twelve vertices.
The binary bound gives respectively $16,16,64$.
For $p=3$, the odd budget leaves central exponents at most
$3,3,4$; exponent three is excluded by
Theorem~\ref{thm:sylow-budgets}(iii). The resulting ternary maxima
are $9,9,81$. For a nonabelian Sylow subgroup $P=F_{(p)}$ with
$p\ge5$, the inequalities
$\omega(P)-1\ge p\log_pz(P)/2$ and $\omega(P)\ge p+1$
leave $25$ at $p=5$, and $49$ at $p=7$ only for the cap eight.
Thus the maxima at $5,6,8$ are $16,25,81$.
The cap seven follows from eight and the ternary example.
The models of Section~\ref{sec:intro} attain every bound.
\end{proof}

\subsection{Binary rigidity}
The numerical profile leaves one extra vertex at even binary
caps. The next theorem rules it out unless the commutator form is
scalar. Its proof uses only the bounds already established.

\begin{theorem}[Uniform rigidity at a binary-square index]
\label{thm:uniform-rigidity}
Let $G$ be finite and suppose $z(G)=2^n$, where $n$ is even
and $n\ge6$. Exactly one of the following holds:
\[
\begin{aligned}
&G'\cong C_2,\quad G/\Z(G)\cong\F_2^n,
\quad \omega(G)=n+1;\\
&\omega(G)\ge n+3.
\end{aligned}
\]
In particular $\omega(G)=n+2$ is impossible.
\end{theorem}
\begin{proof}
Every odd-order element lies in $\Z(G)$, so $G=P\times O$
with $P$ a $2$-group and $O$ central of odd order. Thus we may
replace $G$ by $P$. Suppose $\omega(P)\le n+2$.
We first show that every scalar radical is $\Z(P)$; affine
inflation will then force the derived group to have order two.
For any scalar radical $R$, the interaction bound gives
$dm+h\le1$, hence $m=0$ and $h\le1$.

\emph{Excluding a noncentral radical.}

If $R>\Z(P)$, then $R$ is abelian and $[R,P]=\langle c\rangle$
has order two, hence is central. For $u\in R$, the identity
$[u^2,g]=[u,g]^2=1$ gives $R^2\le\Z(P)$. Thus the pairing
\[
 R/\Z(P)\times V\longrightarrow\langle c\rangle,
 \qquad (u\Z(P),gR)\longmapsto[u,g],\qquad V=P/R,
\]
is bilinear with trivial left radical, so $1\le s\le d$.
The right radical of this pairing is a proper subspace of $V$.
Choose a hyperplane containing it and a scalar clique of size $d+1$
meeting that hyperplane once. Each of the other $d$ vertices acts
nontrivially on $R$, so its fibre contains at least two
noncommuting elements. The remaining fibre contributes one,
giving $\omega(P)\ge2d+1$. If $s\le d-2$, this is at least $n+3$.
Since $n=d+s$ and both $n,d$ are even, the remaining case is
$s=d$, making the pairing perfect and $n=2d$, $d\ge4$.

Choose a symplectic clique $a_1,\ldots,a_{d+1}$ in $V$ and
$v\notin\{0,a_1,\ldots,a_{d+1}\}$.
Let $v_1,\ldots,v_{2d+3}$ consist of two copies of each $a_j$
and one copy of $v$, and lift them to $g_i\in P$.
The perfect pairing realizes every $f_i\in V^*$ by $b_i\in R$ with
$[b_i,g]=c^{f_i(gR)}$. Moreover $[R,P]\le\Z(P)$ implies
$[R,P']=1$: for $u\in R$, the map $g\mapsto[u,g]$ is a
homomorphism to the centre and therefore kills $P'$. Hence
\[
 [g_ib_i,g_jb_j]=[g_i,g_j]c^{f_i(v_j)+f_j(v_i)}.
\]
Nonzero scalar pairings remain nontrivial because
$c\in[R,P]$ and $R$ is the scalar radical, so $\chi(c)=0$.
For every zero scalar pairing whose commutator lies in
$\langle c\rangle$, prescribe its correction exponent to make
the displayed commutator equal to $c$; prescribe zero on the
other zero pairings.
To meet the prescriptions, denote the functionals for the two copies
of $a_j$ by $f_{j,1},f_{j,2}$, and set the functional for $v$ equal
to zero. Distinct doubled families have scalar pairing one, so
there is no constraint between them. Within the $j$th family the
only constraint is
\[
 f_{j,1}(a_j)+f_{j,2}(a_j)=\varepsilon_j,
\]
where $\varepsilon_j\in\F_2$ is the prescribed correction for its
repeated pair. If $b_\chi(a_j,v)=0$, the two additional constraints
are $f_{j,\ell}(v)=\eta_{j,\ell}$ for $\ell=1,2$, where
$\eta_{j,\ell}\in\F_2$ are the prescribed corrections for the
pairs with $v$. If $b_\chi(a_j,v)\ne0$, there is no constraint
on these values. Choose $f_{j,1}(a_j)=0$ and
$f_{j,2}(a_j)=\varepsilon_j$. Since $a_j$ and $v$ are linearly
independent, each choice of the required values on these two
vectors extends to a functional on $V$. The choices for distinct
families are independent, so all prescriptions hold simultaneously.
The corrected lifts therefore form a clique of size
$2d+3=n+3$, again a contradiction.

\emph{Identifying the commutator form.}
Thus every scalar radical is $\Z(P)$. Since it contains $P'P^2$,
the group has class two and elementary abelian central quotient;
$P'$ is also elementary abelian. If $\dim_{\F_2}P'\ge2$,
Lemma~\ref{lem:affine-inflation} gives
$\omega(P)\ge2n-3\ge n+3$. Consequently $P'\cong C_2$,
and its nondegenerate symplectic form gives $\omega(P)=n+1$.
The two alternatives are exclusive, as claimed.
\end{proof}

\begin{remark}[The dimension threshold]
\label{rem:rank-four-exception}
The restriction $n\ge6$ is necessary: $H=\mathrm{UT}_3(\F_4)$ has
$H'=\Z(H)\cong C_2^2$, $z(H)=16$ and $\omega(H)=5$.
Indeed, its commutator is the determinant of the two coordinates
in $\F_4^2$, so a maximum clique has one representative from
each of the five points of $\mathbb P^1(\F_4)$.
Thus the scalar dichotomy fails when $n=4$.
\end{remark}

\section{Centralizer reduction and saturation}\label{sec:descent}
The nilpotent case is settled. To control a general group, we must
bound a centralizer index and the centrality lost outside that
subgroup. We first quantify the loss, then select the centralizer
and obtain the numerical bound used later. Two saturation
principles will exclude the finite cases left by that bound.

\subsection{Charging a loss of centrality}
For $H\le G$ and $\Z(G)\le A\le\Z(H)$, the index $[A:\Z(G)]$ measures how
much of $A$ is central in $H$ but not in $G$. Adjoining generators
from $H$ up to $G$ removes this excess. At each nontrivial step,
Lemma~\ref{lem:twist} forces new clique vertices. No normality of
$A$ in $G$ is required.

For a positive integer $b$, let $v_p(b)$ be the exponent of $p$ in $b$.
Define the additive prime cost and, for integers $D\ge0$, its product envelope by
\[
 \sigma(b)=\sum_p p\,v_p(b),\quad \sigma(1)=0,\qquad
 \mathcal P(D)=\max\{b\in\mathbb Z_{\ge1}:\sigma(b)\le D\}.
\]
Since $\sigma(uv)=\sigma(u)+\sigma(v)$, the envelope satisfies
$\mathcal P(D)\mathcal P(E)\le\mathcal P(D+E)$ for $D,E\ge0$.
Splitting a composite factor never increases its additive cost, since
$a+b\le ab$ for $a,b\ge2$. Thus $\sigma(b)$ is the minimum sum of
nontrivial factors with product $b$. Integer break gives
$\mathcal P(0)=\mathcal P(1)=1$ and, for $D\ge2$, the usual product of
$2$'s and $3$'s: $3^{D/3}$ if $D\equiv0$, $4\cdot3^{(D-4)/3}$ if
$D\equiv1$, and $2\cdot3^{(D-2)/3}$ if $D\equiv2\pmod3$.

\begin{theorem}[Central-loss calculus]\label{thm:budget}
Let $H\le G$, $\Z\le A\le\Z(H)$ and put $b=[A:\Z]$. Then
\begin{equation}\label{eq:central-budget}
 \omega(G)-\omega(H)\ge\sigma(b).
\end{equation}
More explicitly, for $G=\langle H,g_1,\ldots,g_t\rangle$, put
$H_0=H$, $H_i=\langle H,g_1,\ldots,g_i\rangle$,
$A_i=A\cap\Z(H_i)$ for $0\le i\le t$, and
$d_i=[A_{i-1}:A_i]$ for $1\le i\le t$. Then
\begin{equation}\label{eq:sequential}
 \sum_{d_i>1}d_i\le N-\omega(H),\qquad \prod_i d_i=b.
\end{equation}
For distinct primes $p,q\mid b$, one moreover has
\begin{equation}\label{eq:mixed-surcharge}
 N-\omega(H)\ge\sigma(b)+pq-p-q.
\end{equation}
If $b=36$, this lower bound improves to twelve.
No normality of $A$ in $G$ is assumed.
\end{theorem}
\begin{proof}
Since $A_i=C_{A_{i-1}}(g_i)$, twisting adds at least $d_i$
vertices at each nontrivial step. Sum these increases and use
$A_t=\Z$. Splitting the factors $d_i$ into primes proves
\eqref{eq:central-budget}.

For the surcharge choose $g$ outside the union of the two proper
subgroups $C_G(A_{(p)})$ and $C_G(A_{(q)})$.
Then $d=[A:C_A(g)]$ is divisible by $pq$.
Twisting first by $g$ and applying the ordinary cost to the remaining
loss gives $N-\omega(H)\ge d+\sigma(b/d)$.
Splitting $d$ while retaining $pq$ as one factor gives
$d\ge pq+\sigma(d/pq)$, proving \eqref{eq:mixed-surcharge}.
When $b=36$, the first mixed index is $6,12,18$ or $36$.
Only six could cost less than twelve; it leaves a loss six, whose
mixed surcharge costs six more.
\end{proof}

Every central loss with available budget $D$ satisfies
\begin{equation}\label{eq:coarse-loss}
 b\le\mathcal P(D)\le3^{D/3}.
\end{equation}
We set $\mathcal P(D)=0$ for $D<0$.
At budgets $D=5,7,11$, the mixed surcharge sharpens this to
$b\le5,9,36$, respectively: mixed losses have
$\sigma(b)\le D-1$, while the largest pure prime powers of
cost at most $D$ are $5,9,32$.

For later use, put $K=C_G(L)$, $a=[G:K]$, $b=[\Z(K):\Z]$
and $q=\omega(L/\Z(L))$. Then
\begin{equation}\label{eq:centralizer-accounting}
 z(G)=a\,z(K)b,\qquad N\ge q\omega(K)+\sigma(b).
\end{equation}
Indeed $\Z(KL)=\Z(K)$: a central element of $KL$ centralizes
$L$ and hence belongs to $K$. Apply \eqref{eq:product-amplification}
and the central-loss bound with $H=KL$, $A=\Z(K)$.

\subsection{Fourfold centralizer selection}
\begin{lemma}[Complement rigidity]\label{lem:complement}
If $k\ge1$ is an integer, $D<G$ and $\omega(G\setminus D)\le k$,
then either $C_G(x)$ is abelian of index at most $k$ for some
$x\notin D$, or some $x\notin D$ satisfies $[G:C_G(x)]\le k-1$.
\end{lemma}
\begin{proof}
Suppose all those indices are at least $k$. Put $S=G\setminus D$
and $d=[G:D]$. For $C=C_G(x)$, $x\in S$, one has
$q=[C:C\cap D]\le d$ and hence
\[
 |C\cap S|=|C|(1-1/q)\le |S|/k.
\]
A maximum clique in $S$ covers $S$ by at most $k$ such sets.
Thus exactly $k$ sets partition $S$, each with centralizer of
index $k$. The elements of each part commute pairwise: otherwise two of its vertices
could replace its selected vertex, enlarging the clique.
Moreover $C\setminus D$ generates $C$, since for $s\in C\setminus D$
and $a\in C\cap D$, both $s,sa$ belong to that part. Thus $C$ is abelian.
\end{proof}

\begin{proposition}[Fourfold descent]\label{prop:descent}
If $G$ is nonnilpotent, one of the following holds:
\begin{enumerate}[label=\textup{(\Alph*)}]
\item There is $x\notin F(G)$ such that $C_G(x)$ is abelian
 and $[G:C_G(x)]\le N$.
\item There are conjugate elements $x,y$, with $x\notin F(G)$,
 such that $L=\langle x,y\rangle$ is nonnilpotent. For
 $K=C_G(L)$, $r=|x^G|$ and $a=[G:K]$,
 \begin{equation}\label{eq:fourfold-descent}
 r\le N-1,\qquad r\mid a\le r(r-1),\qquad
 \omega(KL)\ge4\omega(K).
 \end{equation}
\end{enumerate}
Here $F(G)$ denotes the Fitting subgroup.
\end{proposition}
\begin{proof}
Apply Lemma~\ref{lem:complement} to $D=F(G)$.
Outside its abelian alternative, choose $x\notin F(G)$ with
$r\le N-1$. Baer--Suzuki in its nilpotent formulation
\cite[p.~281]{GuralnickRobinson1993} supplies a conjugate $y$
for which $L=\langle x,y\rangle$ is nonnilpotent.
The group $C_G(x)$ fixes $x$ in $x^G$, so its orbit on $y$
has size at most $r-1$. Its stabilizer is $K$, proving the
index bounds. Finally $L/\Z(L)$ is nonnilpotent and has a
four-point clique by Lemma~\ref{lem:small-quotients}.
Apply \eqref{eq:product-amplification}.
\end{proof}

\subsection{Losses for an abelian centralizer}
\begin{lemma}[Abelian centralizer losses]\label{lem:marked-loss}
Let $L$ be nonnilpotent, $K=C_G(L)$ abelian, and put
$C=C_G(K)$, $a=[G:K]$, $b=[K:\Z]$ and $h=\omega(C)$.
Then $\Z(C)=K$, $C/K$ is nonnilpotent, and
\begin{equation}\label{eq:double-centralizer-profile}
 h\ge4,\qquad h+\sigma(b)\le N.
\end{equation}
If $h=4$, then $C/K\cong S_3$, $6\mid a$, and $C$ has an
abelian subgroup of loss $3b$.
Consequently, for $N\le k$ and $k\ge5$,
\begin{equation}\label{eq:abelian-marked}
 b\le J(k):=\max\{\mathcal P(k-5),3^{\lfloor(k-4)/3\rfloor}\};
 \qquad N\le12\Longrightarrow b\le9.
\end{equation}
\end{lemma}
\begin{proof}
Since $L\le C$, a central element of $C$ lies in $C_G(L)=K$;
the reverse inclusion follows from the definition of $C$.
Also $C/K$ contains the nonnilpotent subgroup $L/\Z(L)$,
so $h\ge4$. Apply the central-loss budget to $K\le\Z(C)$.
If $h=4$, Lemma~\ref{lem:small-quotients} leaves only $S_3$.
The preimage of its subgroup of order three is abelian, since
its quotient by the central subgroup $K$ is cyclic.

For $h\ge5$, $b\le\mathcal P(k-5)$.
For $h=4$, the abelian loss $3b$ has budget $k-1$.
If $b$ is a power of three, its ordinary cost gives the second
term in \eqref{eq:abelian-marked}. Otherwise the mixed surcharge
gives $\sigma(b)+4\le k-1$, and hence the first term.
At cap twelve the budget-seven refinement gives $b\le9$ when
$h\ge5$. For $h=4$, only $b=10,11,12$ remain to exclude:
the losses $30,33,36$ cost at least $17,33,12$, exceeding eleven.

\end{proof}

\subsection{Near-minimal ternary losses}\label{sec:orbit-collapse}
The following theorem uses the prime-index result of Bergman and
Lenstra. A noncentral commutator subgroup of order three would
produce a local twelve-point clique; the central cost of the
remaining dimensions then exceeds the available budget.
This is the divisibility input for the ternary losses left by
Lemma~\ref{lem:marked-loss}.

\begin{theorem}[Near-minimal ternary losses]\label{thm:ternary-saturation}
Let $\Z\le A\le G$ be abelian, with
$[A:\Z]=2^\varepsilon3^r$, $r\ge2$ and $\varepsilon\in\{0,1\}$.
If $N\le3r+3+2\varepsilon$, then $3^r\mid[G:A]$.
\end{theorem}
\begin{proof}
\emph{First force orbits of size at most three.}
Put $D=\Z A_{(3)}$. If $[D:C_D(g)]=3^e\ge9$, taking $g$
first in the central-loss sequence for $A$ and splitting off
any factor two costs at least
\[
 3^e+3(r-e)+2\varepsilon\ge3r+3+2\varepsilon\ge N,
\]
exceeding its budget $N-1$; the splitting uses
$2\cdot3^e\ge3^e+2$ if the binary part also moves.
Thus every $D$-orbit in $G$ has size one or three.

\emph{Apply Bergman--Lenstra.}
Put $E=D/\Z$ and $T=G\rtimes E$, where $d\Z\in E$ acts
on $G$ by the inner automorphism induced by $d$. This action is
well defined and faithful, since its kernel is $C_D(G)/\Z=1$.
We identify $G$ and $E$ with the normal factor and the complement
in $T$, respectively. If a normal subgroup $V$ of $T$ is contained
in $E$, then $[V,G]\le V\cap G=1$; faithfulness gives $V=1$.
Thus $E$ has trivial core in $T$.
For $g\in G$, the orbit bound above gives
$[E:E\cap E^g]=[D:C_D(g)]\in\{1,3\}$.
Since $E$ normalizes itself and $T=GE$, these are also all the
possible indices for conjugation by $T$. Both values occur:
otherwise the action of the nontrivial group $E$ on $G$ would
be trivial.
By \cite[Theorem~2]{BergmanLenstra1989}, a normal subgroup $M$ of $T$
either lies in $E$ with $E/M$ a transitive permutation group of
degree three, or contains $E$ with $[M:E]=3$.
The first case forces $M=1$ and $|E|=3$, since a faithful abelian
transitive group is regular, contrary to $r\ge2$. Thus the second holds, and
$U=M\cap G$ is normal of order three, since $M=UE$ and
$[M:E]=3$. As $E\le M$ and both $M$ and $G$ are normal
in $T$, we have $[G,E]\le M\cap G=U$. This commutator
subgroup is nontrivial by faithfulness, so it is all of $U$.
The action defining the semidirect product identifies it with
$[D,G]$ inside $G$. Finally $E$ centralizes $U$, because a
$3$-group acts trivially on $C_3$, whose automorphism group has
order two. The same holds for conjugation by $D$.

\emph{A local twelve-point clique forces centrality.}
We claim that $U\le\Z$. Otherwise $C_G(U)$ has index two;
choose $g\notin C_G(U)$. Since $D$ is abelian and centralizes
$U$, the group $B=DU$ is abelian. Moreover $[D,G]=U$ and
normality of $U$ give $[B,G]=U\le B$, so $B$ is normal.
Also $[B:\Z]=3^t$ for $t\in\{r,r+1\}$.
There is $u\in C_G(U)$ not centralizing $B$: otherwise
$C_B(g)=C_B(G)=\Z$, whereas $[B:C_B(g)]=3$ and $t\ge2$.
The index-three subgroups $C_B(u),C_B(g)$ intersect in $B_0$
with $[B:B_0]=9$, because only $C_B(u)$ contains $U$.
Put $H_0=\langle B,u\rangle$, $H=\langle H_0,g\rangle$ and
$C=C_{H_0}(B)$. The action of $H_0$ on $B$ is generated by
the order-three action of $u$, so $[H_0:C]=3$; also $B_0\le\Z(H)$.

For $w\in gH_0$, we have $C_{H_0}(w)\le C$.
Indeed, if $h\in H_0$ commutes with $w$, then for every $b\in B$,
\[
 [b,h]=[b^w,h]=[b,h]^w=[b,h]^{-1};
\]
the first equality uses $b^w\in bU$ and $[U,H_0]=1$, and
$w$ inverts $U$. Thus $[b,h]=1$.
Also $[B:C_B(w)]=3$, so
\[
 [H_0:C_{H_0}(w)]\ge[H_0:C][B:C_B(w)]=9.
\]
Since this index bound holds for every $w\in gH_0$, the averaging
assertion in Lemma~\ref{lem:coset-calculus} gives a nine-point
clique in $gH_0$.
Choose $a\in B$ with $[u,a]\ne1$. The three elements $ua^j$,
$0\le j<3$, form a clique in $H_0\setminus C$: the nontrivial
commutator $[u,a]\in U$ is central in $H_0$ and has order three.
They are completely adjacent to $gH_0$, because
$C_{H_0}(w)\le C$ for every $w\in gH_0$. Hence $\omega(H)\ge12$.
Applying the central-loss bound to $B_0\le\Z(H)$ now gives
\[
 N\ge12+\sigma([B_0:\Z])
   =12+3(t-2)\ge3r+6>N,
\]
a contradiction. Therefore $U\le\Z$, so $D$ is normal and $[D^3,G]=1$ makes
$D/\Z\cong\F_3^r$.

\emph{Read off the divisibility.}
The commutator map is the homomorphism
\[
 G\longrightarrow\operatorname{Hom}(D/\Z,U),
 \qquad g\longmapsto(d\Z\longmapsto[d,g]).
\]
Identify $U$ with the additive group of $\F_3$. The image is
then a subspace of the dual of $D/\Z\cong\F_3^r$, with zero
annihilator because $C_D(G)=\Z$. A proper subspace of this
dual has a nonzero annihilator, by the dimension formula for
the evaluation pairing. Thus the image is the whole dual.
Its kernel is $C_G(D)$, so
$[G:C_G(D)]=3^r$. Since $D\le C_G(D)$, we obtain
$3^r\mid[G:D]=2^\varepsilon[G:A]$, proving the assertion.
\end{proof}

\paragraph{Lifting the ternary loss.}
Return to the notation of Lemma~\ref{lem:marked-loss}, and put $m=[C:K]$.
If $b=2^\varepsilon3^s$, $\varepsilon\in\{0,1\}$, $s\ge1$,
$3\mid m$ and $N\le3s+6+2\varepsilon$, then
\begin{equation}\label{eq:ternary-lift}
 3^{s+2}\mid a.
\end{equation}

For \eqref{eq:ternary-lift}, lift an element of order three in
$C/K$ to $y$. Then $A=K\langle y\rangle$ is abelian with loss
$3b$ and index $a/3$ in $G$. Theorem~\ref{thm:ternary-saturation}
gives $3^{s+1}\mid a/3$, as required.

\subsection{Saturation of an affine conjugacy class}
In this subsection we write the conjugation action on a class
using left conjugation. Choose a bijection $u\mapsto x_u$ from a finite abelian group
$A$ onto the class $X$.
The law is affine if
\[
 x_u x_v x_u^{-1}=x_{(1-f)u+fv},
 \qquad f\in\operatorname{Aut}(A),\quad 1-f\text{ invertible}.
\]
Thus $u,v$ are labels in $A$, whereas $x_u,x_v$ are elements of $G$.
Such a class is connected: the subgroup it generates acts transitively on it.

\begin{lemma}[Affine saturation]\label{lem:affine-saturation}
Suppose a conjugacy class $X$ of $G$ has the affine law above,
$|A|=q>1$, and put
$R=\langle X\rangle$. Then $\omega(G)\ge q+1$.
If $\omega(G)\le q+1$, then
\[
 C_G(R)=\Z(G),\qquad z(G)\mid q|\operatorname{Aut}(A)|.
\]
\end{lemma}
\begin{proof}
Conjugation identifies $R/\Z(R)$ with the inner group of $X$.
For the connected affine law, the standard descriptions give
$R/\Z(R)\cong A\rtimes\langle f\rangle$ and
$\operatorname{Aut}(X)=A\rtimes C_{\operatorname{Aut}(A)}(f)$;
see \cite[Lemma~1.9, Corollary~1.25]{AndruskiewitschGrana2003}.
The $q$ elements $(a,f)$, $a\in A$, are pairwise noncommuting,
and every $(w,1)$ with $w\ne0$ fails to commute with all of them,
since $1-f$ is invertible. Thus $\omega(R/\Z(R))\ge q+1$.

Put $E=C_G(R)$. If $N\le q+1$, \eqref{eq:centralizer-accounting}
gives
\[
 q+1\ge N\ge(q+1)\omega(E)+\sigma([\Z(E):\Z]).
\]
Hence $E=\Z$. Conjugation embeds $G/\Z$ in
$\operatorname{Aut}(X)$, so
$z(G)\mid q|C_{\operatorname{Aut}(A)}(f)|$,
which implies the required divisibility.
\end{proof}

\begin{corollary}[Small conjugacy classes]\label{cor:small-classes}
Let $x\notin F(G)$ and $r=|x^G|$.
If $r=p$ is prime, then $N\ge p+1$; if also $N\le p+1$,
then $z(G)\mid p(p-1)$.
If $r=9$ and $N\le10$, either there is a conjugate $y$ with
$\langle x,y\rangle$ nonnilpotent and $[G:C_G(x,y)]\le18$,
or $z(G)\mid54$ or $z(G)\mid432$.
\end{corollary}
\begin{proof}
Set $X=x^G$ and $R=\langle X\rangle$. Normality of $R$ makes
all its orbits on $X$ equally large, so $|x^R|\mid r$.
This size is not one: otherwise $X\subseteq\Z(R)$, so $R$ is
abelian and $x\in F(G)$.
If $r=9$ and $|x^R|=3$, then $F(R)\le F(G)$ by characteristicity.
Baer--Suzuki in $R$ supplies $y\in x^R$ with $\langle x,y\rangle$
nonnilpotent. Since $C_G(x)$ preserves $x^R$ and fixes $x$,
$[G:C_G(x,y)]=9|y^{C_G(x)}|\le18$.

In the remaining cases, $X=x^R$ is connected. The prime-order theorem of
Etingof--Soloviev--Guralnick \cite{EtingofGuralnickSoloviev2001}
and Gra\~na's prime-square theorem \cite[Proposition~3.10]{Grana2004}
make its law affine. Apply Lemma~\ref{lem:affine-saturation}.
For prime $p$, $A=C_p$ has $|\operatorname{Aut}(A)|=p-1$;
for nine, $A=C_9$ or $C_3^2$ has automorphism group of order six
or $48$, giving the two stated divisibilities.
\end{proof}
\section{Completion of the exact profile}\label{sec:completion}
We first establish the small caps and isolate the finite exceptions.
Those exceptions are grouped by the mechanism that excludes them:
ternary divisibility at caps eight, fourteen and sixteen, and affine
saturation at cap ten. After this finite block, one numerical envelope
completes the proof by induction.

\subsection{The global upper bound}
\begin{proof}[Proof of Theorem~\ref{thm:global}]
The models in Section~\ref{sec:intro} give every lower bound.
Proposition~\ref{prop:nilpotent-profile} gives every nilpotent upper
bound, and Lemma~\ref{lem:small-quotients} settles all caps $k\le4$. We therefore consider nonnilpotent groups. It is enough to
prove the bounds at cap five and at every even cap $k\ge6$:
each odd cap $2j-1\ge7$ then follows from cap $2j$, since
$M(2j-1)=M(2j)$.

\medskip\noindent
\textbf{Step 1. A common budget for the two cases.}
In case \textup{(A)} of Proposition~\ref{prop:descent}, put
$A=C_G(x)$, $a=[G:A]$ and $b=[A:\Z]$. Then
\[
 a\le k,\qquad \sigma(b)\le k-1,\qquad z(G)=ab.
\]
In case \textup{(B)}, retain the selected class size $r$ and put
$a=[G:K]$, $t=\omega(K)$ and $b=[\Z(K):\Z]$. The accounting
identity and fourfold amplification give
\begin{equation}\label{eq:profile-accounting}
 z(G)=a\,z(K)b,\qquad 4t+\sigma(b)\le k,\qquad
 r\mid a\le r(r-1),\quad r\le k-1.
\end{equation}
If $K$ is abelian, Lemma~\ref{lem:marked-loss} gives $b\le J(k)$.
If $K$ is nonabelian, then $t\ge3$. Whenever the profile is already
known at $t$, it gives $z(K)\le M(t)$.

Here is the numerical fact that lets us absorb this internal index.
For $t=3,\ldots,8$, the respective values of $\sigma(M(t))$ are
$4,6,8,10,12,12$; for $t\ge9$, the binary formula gives
$\sigma(M(t))=4\lfloor(t-1)/2\rfloor$. Thus
$\sigma(M(t))\le2(t-1)$ for every $t\ge3$, and
\begin{equation}\label{eq:cost-absorption}
 \begin{aligned}
 z(K)b
 &\le M(t)\mathcal P(k-4t)\\
 &\le\mathcal P(k-2t-2)
 \le\mathcal P(k-8)\le J(k).
 \end{aligned}
\end{equation}
The second line uses the additive prime cost of the numerical
majorant in the first line. Consequently, once the bound at $t$ is
available, every selected centralizer satisfies $z(G)\le aJ(k)$.
For all caps $k\le16$, we have $t\le4$, so the required
bound at $t$ is already supplied by
Lemma~\ref{lem:small-quotients}. For the later caps it will be
supplied by induction.

\medskip\noindent
\textbf{Step 2. The base caps five and six.}
At cap five, the cases $N\le4$ are already settled. If $N=5$,
Lemma~\ref{lem:coset-calculus} gives an irredundant cover by five
centralizers with intersection $\Z$, so $z(G)\le16$ by
\cite[Theorem~1.2]{BryceFedriSerena1997}.

At cap six, case \textup{(B)} has $4t\le6$, so $K$ is abelian.
Since $J(6)=1$, this gives $z(G)\le20<25$.
In case \textup{(A)}, $a\le6$ and $b\le5$, so an excess over
$25$ requires $(a,b)=(6,5)$.

For $(a,b)=(6,5)$, the quotient $Q=G/\Z$ has order thirty and
a normal Sylow $5$-subgroup. Indeed, otherwise its six Sylow $5$-subgroups
supply twenty-four nonidentity elements, forcing the Sylow $3$-subgroup
to be normal. The quotient of order ten then gives a normal cyclic
subgroup of order fifteen; its Sylow $5$-subgroup is characteristic and
therefore normal in $Q$, a contradiction.
Hence $A/\Z\cong C_5$ is normal. Let $H$ be the kernel of the
action of $G$ on $A/\Z$. Since $A=C_G(x)$, the map
\[
 H/A\longrightarrow\Z,\qquad hA\longmapsto[x,h]
\]
is injective, and its image has exponent five because $x^5\in\Z$.
Its order also divides six, so $H=A$. This embeds the group
$G/A$, of order six, in $\operatorname{Aut}(C_5)$, of order four.
Thus $\zeta(5)=M(5)=16$ and $\zeta(6)=M(6)=25$.

\medskip\noindent
\textbf{Step 3. The finite caps eight through sixteen.}
The following table records upper bounds for $z(G)$ in each
alternative. An entry at most $M(k)$ closes that branch; a larger
entry indicates a case still to be examined. A dash means that
$4\omega(K)\le k$ excludes a nonabelian $K$.
For case \textup{(A)}, use $k\mathcal P(k-1)$,
sharpened at budgets seven and eleven to losses nine and
thirty-six. For nonabelian $K$, use
$(k-1)(k-2)\mathcal P(k-8)$ from
\eqref{eq:cost-absorption}. For abelian $K$, use
$(k-1)(k-2)J(k)$, with the sharper $b\le9$ at cap twelve
from Lemma~\ref{lem:marked-loss}.

\begin{center}
\begin{minipage}{\linewidth}
\centering
\small
\refstepcounter{table}\label{tab:finite-cap-screen}
\begin{tabular}{c|r|r|r|r}
 $k$&$M(k)$&case \textup{(A)}
 &\shortstack{case \textup{(B)}\\$K$ nonabelian}
 &\shortstack{case \textup{(B)}\\$K$ abelian}\\\hline
 8&81&72&---&126\\
 10&256&270&---&648\\
 12&1024&432&440&990\\
 14&4096&1512&1404&4212\\
 16&16384&3888&3780&17010
\end{tabular}
\par\smallskip
\textbf{Table~\thetable.} Preliminary bounds for the finite caps. The three right-hand
columns distinguish case \textup{(A)} and the two possibilities
for $K$ in case \textup{(B)}; the entries are bounds, not attained values.
\end{minipage}
\end{center}

The case $k=12$ is therefore settled. In case \textup{(A)} the only possible excess is at cap ten, where
$a\le10$, $b\le27$ and $ab>256$ require $a=10$ and $b=26$ or
$27$. The loss $26$ costs fifteen, while loss $27$ is excluded by
Theorem~\ref{thm:ternary-saturation}, which would give $27\mid10$.
All remaining excesses thus have abelian selected $K$.

We now treat caps eight, fourteen and sixteen together, and then
return to ten. This order has no inductive dependency: throughout
this finite block the only profile bounds used are those at caps
at most five, already proved above.

\medskip\noindent
\textbf{Step 4. Three cases excluded by ternary divisibility.}
Assume $z(G)>M(k)$ in case \textup{(B)}, with $K$ abelian.
At cap eight, if $r=7$, Corollary~\ref{cor:small-classes} gives
$z(G)\mid42$. Otherwise $r\le6$. Since $b\le3$ and
$r\mid a\le r(r-1)$, the inequality $ab>81$ then forces
$(a,b)=(30,3)$.

At cap fourteen, $a\le156$, $b\le27$ and $ab>4096$ force
$b=27$ and $a\ge152$. Hence $r=13$ and $a=156$.
At cap sixteen, $a\le210$, $b\le81$ and $ab>16384$ leave
$79\le b\le81$. But $\sigma(b)\le12$ excludes $79$ and $80$,
so $b=81$; the index constraints then give $r=15$ and $a=210$.

Put $C=C_G(K)$ and $m=[C:K]$. At caps eight and fourteen,
Lemma~\ref{lem:marked-loss} gives
$\omega(C)\le k-\sigma(b)=5$. Since $\Z(C)=K$, the cap-five
bound from Step~2 yields $m\le16$. Also $m\mid a$ and $C/K$ is
nonnilpotent. For $a=30$ this leaves $m=6$ or $10$; for $a=156$
it leaves $m=6$ or $12$. The order-ten possibility is dihedral
and has a six-point clique (five reflections and a nontrivial
rotation), contradicting $\omega(C)\le5$.
At cap sixteen, $4\le\omega(C)\le16-12=4$, so
$C/K\cong S_3$ and $m=6$.

In every remaining row, $3\mid m$. The ternary lifting criterion
\eqref{eq:ternary-lift} now gives the contradictions in
Table~\ref{tab:ternary-contradictions}. The criterion requires
respectively $8\le9$, $14\le15$ and $16\le18$.
\begin{center}
\begin{minipage}{\linewidth}
\centering
\refstepcounter{table}\label{tab:ternary-contradictions}
\begin{tabular}{c|r|r|c|c}
 $k$&$a$&$b$&$m$&required divisibility\\\hline
 8&30&3&6&$27\mid30$\\
 14&156&27&$6$ or $12$&$243\mid156$\\
 16&210&81&6&$729\mid210$
\end{tabular}
\par\smallskip
{\small\textbf{Table~\thetable.} Divisibility obstructions in the three ternary cases.
Each divisibility is forced by \eqref{eq:ternary-lift} and is
incompatible with the displayed value of $a$.}
\end{minipage}
\end{center}

\medskip\noindent
\textbf{Step 5. The remaining cases at cap ten.}
It remains to consider cap ten with $K$ abelian. Here $a\le72$,
$b\le9$ and $ab>256$, so $b\ge4$. Again put $C=C_G(K)$.
Loss five is impossible: lift an element of prime order $q$ in
$C/K$. Together with $K$, its lift generates an abelian subgroup
of loss $5q$. The central-loss calculus assigns this loss a cost
of at least ten, exceeding the available budget nine.

For $b=6,7,8,9$, the central-loss cost is at least six; for $b=6$
this uses the mixed surcharge, and for the other values the ordinary
cost suffices. Thus $\omega(C)=4$, unless loss seven already
contradicts $4+\sigma(b)\le10$. The small-quotient lemma gives
$C/K\cong S_3$ and $6\mid a$. Lifting its subgroup of order three
gives an abelian subgroup of loss $3b$. If $b=8$, this is loss
$24$, whose mixed cost is ten, again impossible. If $b=9$,
\eqref{eq:ternary-lift} gives $81\mid a$, contrary to $a\le72$.
If $b=6$, the same criterion gives $27\mid a$; together with
$6\mid a$, this forces $a=54$, and the class-index constraints
force $r=9$.

Finally, if $b=4$, excess requires $a\ge65$. The constraints
$r\le9$, $r\mid a\le r(r-1)$ then give $r=9$ and $a=72$.
We have reduced the entire cap-ten problem to
\[
 (r,a,b,z(G))=(9,54,6,324)
 \quad\hbox{or}\quad(9,72,4,288).
\]
Apply Corollary~\ref{cor:small-classes}. Either it supplies a new
selected centralizer of index at most eighteen, giving
$z(G)\le18J(10)=162$, or it gives $z(G)\mid54$ or $z(G)\mid432$.
Both values $324$ and $288$ exceed $162$, and neither divides
$54$ or $432$. This settles cap ten and completes every cap
through sixteen.

\medskip\noindent
\textbf{Step 6. Induction for all larger even caps.}
Proceed by strong induction on the even caps $k\ge18$.
In case \textup{(B)}, $t\le k/4<k$;
the profile bound needed in Step~1 follows from smaller caps.
If $t$ is odd and $t\ge7$, use the already proved even cap
$t+1\le k-2$; otherwise the bound at $t$ is already available. Since $J(k)\le3^{(k-4)/3}$, we obtain
\[
 z(G)\le T(k):=(k-1)(k-2)3^{(k-4)/3}.
\]
The same envelope covers case \textup{(A)}:
$z(G)\le k\mathcal P(k-1)\le k3^{(k-1)/3}$, and
$3k\le(k-1)(k-2)$.

At the first cap,
\[
 T(18)=272\cdot81\cdot3^{2/3}
       <272\cdot81\cdot\frac94=49572<65536=L_{18}.
\]
The rational factor in the following ratio decreases with $k$, so
for every even $k\ge18$,
\[
 \frac{T(k+2)}{T(k)}
 =\frac{k(k+1)}{(k-1)(k-2)}3^{2/3}
 <\frac{171}{136}\frac94<3<4.
\]
Since $L_{k+2}=4L_k$ at the even caps, the envelope remains
strictly below the target. This completes the induction and the
proof of the exact profile.
\end{proof}

\subsection{Extremizers and the sharp constant}
The cap five requires one elementary classification of alternating maps.
It accounts for the two different extremal models at central index sixteen.

\begin{lemma}[Nonsingular alternating maps in dimension four]
\label{lem:four-dimensional-maps}
Let $V=\F_2^4$, let $W\ne0$ be a finite-dimensional $\F_2$-space,
and let $\beta:V\times V\to W$ be alternating and bilinear, with
values spanning $W$. Suppose that $\lambda\circ\beta$ is nondegenerate
for every nonzero $\lambda\in W^*$. Up to linear identifications of
$V$ and $W$, precisely two possibilities occur:
\begin{enumerate}[label=\textup{(\alph*)}]
\item $W=\F_2$ and $\beta$ is a nondegenerate symplectic form;
\item $V=\F_4^2$, $W=\F_4$, and $\beta$ is the determinant over $\F_4$.
\end{enumerate}
\end{lemma}
\begin{proof}
The scalar forms constitute a subspace $S$ of the alternating forms
on $V$, of dimension $\dim W$, since the values of $\beta$ span $W$.
For an alternating matrix $C=(c_{ij})$ of size four, its Pfaffian is
\[
 Q(C)=c_{12}c_{34}+c_{13}c_{24}+c_{14}c_{23},
 \qquad \det C=Q(C)^2.
\]
Consequently $Q$ takes the value one on every nonzero element of $S$.
If $S$ contained a three-dimensional subspace, the sum of $Q$ over
that subspace would be one in $\F_2$, since there are seven nonzero
vectors. But a polynomial of degree at most two on $\F_2^3$ sums
to zero: each of its monomials is one at an even number of vectors.
Thus $\dim W\le2$. Dimension one gives \textup{(a)}.

In dimension two, write $\beta=(b,c)$, where $b,c,b+c$ are
nondegenerate alternating forms, and define $T$ by
$c(x,y)=b(x,Ty)$. Choose a symplectic basis for $b$, so its matrix is
$J=\left(\begin{smallmatrix}0&I_2\\I_2&0\end{smallmatrix}\right)$.
Write the matrix of $c$ as
\[
 C=\begin{pmatrix}aJ_2&B\\B^{\mathsf t}&dJ_2\end{pmatrix},
 \qquad J_2=\begin{pmatrix}0&1\\1&0\end{pmatrix}.
\]
Since $T=JC$, the two-by-two Cayley--Hamilton identity gives
\[
 T^2+\operatorname{tr}(B)T+(\det B+ad)I_4=0.
\]
Indeed, the diagonal blocks follow from Cayley--Hamilton for $B$,
and the off-diagonal blocks use
$B^{\mathsf t}J_2+J_2B=\operatorname{tr}(B)J_2$.
Both $T$ and $T+I_4$ are invertible, since $c$ and $b+c$ are
nondegenerate. Of the four monic quadratics over $\F_2$, only
$t^2+t+1$ can annihilate such an operator. Hence $T^2+T+I_4=0$.

Let $\F_4=\F_2(\tau)$, where $\tau^2+\tau+1=0$, and let $\tau$
act on $V$ as $T$. The operator $T$ is self-adjoint for $b$, because
$c$ is symmetric. The $\F_4$-valued form
$\widetilde\beta=b+\tau c$ is alternating and satisfies
\[
 \widetilde\beta(x,Ty)
 =c(x,y)+\tau\bigl(b(x,y)+c(x,y)\bigr)
 =\tau\widetilde\beta(x,y).
\]
It is therefore bilinear over $\F_4$ and nondegenerate on the
two-dimensional $\F_4$-space $V$. A choice of basis makes it the
determinant, proving \textup{(b)}. Conversely, both listed maps
have the required property. For the determinant, if $v\ne0$, the
map $u\mapsto\det(v,u)$ is onto $\F_4$, so its composition with
any nonzero $\F_2$-linear functional is nonzero. Thus no nonzero
$v$ belongs to the radical of any nonzero scalarization, which
proves the required nondegeneracy.
\end{proof}

\begin{proof}[Proof of Theorem~\ref{thm:equality}]
For $k=1,2$, the absence of a clique of size three forces $G$ to
be abelian. Now assume $k\ge3$.
Every nontrivial value $M(k)$ is a prime power. Hence a maximizer
is nilpotent, with only one noncentral Sylow subgroup.
For $k\in\{3,4,6\}$ and $p=k-1$, the noncyclic central quotient
of order $p^2$ is $C_p^2$. Thus $G$ has class two and two lifts $x,y$ of a basis, with
the centre, generate $G$. Their commutator is nontrivial and satisfies
$[x,y]^p=[x^p,y]=1$, so it generates $G'\cong C_p$.
Conversely, such a scalar symplectic plane gives $z(G)=p^2$ and
$\omega(G)=p+1=k$.

At $k=5$, write $G=P\times A_0$, where $P$ is its Sylow
$2$-subgroup and $A_0$ is central of odd order. Then $z(P)=16$
and $\omega(P)\le5$. The binary budget gives $\omega(P)=5$ and
forces $P$ to have class two with $P/\Z(P)\cong\F_2^4$.
The group $P'$ is elementary abelian, since
$[x,y]^2=[x^2,y]=1$, so the commutator is an alternating bilinear
map $\beta:P/\Z(P)\times P/\Z(P)\to P'$ whose values span $P'$.
For every nonzero scalarization, equality in \eqref{eq:binary-excess}
forces $m=h=0$, hence $R=\Z(P)$. Thus every scalar form of $\beta$
is nondegenerate. Lemma~\ref{lem:four-dimensional-maps} identifies
$\beta$ with either the scalar symplectic form or the determinant
on $\F_4^2$. These are exactly the commutator maps of an
extraspecial group of order $32$ and of $\mathrm{UT}_3(\F_4)$.
Both models have $z=16$ and $\omega=5$, as shown in the introduction.
Their derived subgroups have orders two and four, so they are not
isoclinic. This proves exhaustivity at cap five.

At $k=7,8$, apply Theorem~\ref{thm:sylow-budgets}(iii) to the
Sylow $3$-subgroup. It gives $G'\cong C_3$,
$G/\Z(G)\cong\F_3^4$ and $\omega(G)=7$.
The converse is the scalar $W(3,3)$ calculation.

Finally let $k\ge9$, $r=\lfloor(k-1)/2\rfloor$ and
$z(G)=L_k=2^{2r}$. Theorem~\ref{thm:uniform-rigidity} applies
with $n=2r\ge8$. Its second alternative is impossible because
$2r+3>k$. Hence $G'\cong C_2$ and
$G/\Z(G)\cong\F_2^{2r}$, with $\omega(G)=2r+1$.
The converse follows from the binary symplectic clique calculation.
In each scalar case a symplectic basis identifies the commutator map
with that of the corresponding extraspecial group. Together with
the cap-five classification, this proves all the asserted isoclinisms.
\end{proof}

\begin{proof}[Proof of Corollary~\ref{cor:universal}]
For $N=\omega(G)\ge9$, the exact profile gives
$z(G)\le L_N\le2^{N-1}<(81/64)2^{N-1}$.
For $N\le8$, comparison of the eight values of $M(N)/2^{N-1}$
shows that the maximum is $81/64$, reached only at $N=7$ and
$z(G)=81$. Theorem~\ref{thm:equality}(iv) gives its precise
isoclinism class. The same exact formula yields the limit in
\eqref{eq:exponential-consequences}; the extraspecial binary
family proves optimality of the asymptotic base two.
\end{proof}

\noindent\begin{minipage}{\linewidth}
\section*{Declarations}
\noindent\textbf{Funding.} The author received no funding for this work.
\par\smallskip
\noindent\textbf{Competing interests.} The author declares no competing interests.
\par\smallskip
\noindent\textbf{Data availability.} No empirical datasets were generated or
analysed during this theoretical study.
\par\smallskip
\noindent\textbf{Use of AI-assisted tools.} ChatGPT (OpenAI) was used as an
assistant in checking mathematical proofs and in editing and formatting the
\LaTeX{} manuscript. The author takes full responsibility for the mathematical
content and the final manuscript.
\end{minipage}\par

\end{document}